\documentclass[11pt,reqno]{amsart}
\usepackage[notref,notcite]{}
\usepackage{graphicx}
\usepackage{amssymb}
\usepackage{amsmath}

\begin{document}

\def\?{
?\vadjust{\vbox to 0pt{\vss\hbox{\kern\hsize\kern1em\large\bf ?!}}}}

\title[]{On a holonomy flag of non-holonomic distributions.}

\author{E.~G.~Malkovich}
\address{Sobolev Institute of Mathematics, Russia}
\address{Novosibirsk State University, Russia}
\email{malkovich@math.nsc.ru}

\thanks{The author is supported by
a~Grant of the Russian Federation for the State Support of Researches
(Contract No.~14.B25.31.0029).}


\begin{abstract}

We give definition of a holonomy flag in subRiemannian geometry --- a generalization of a Riemannian holonomy algebra --- and calculate it for the 3D subRiemannian Lie groups. We rewrite and give new interpretation for the Codazzi equations for the $(2,3)$-distributions on the $SU(2)$ and the Heisenberg group.



\end{abstract}

\maketitle

\sloppy

\section{Ideology and motivation.}
In this article we give a new definition of the holonomy flag for a non-holonomic distributions on a subRiemannian manifolds. It generalizes the holonomy algebra for Riemannian manifolds. The usual way to define the holonomy group of Riemannian manifold uses the Levi-Civita connection $\nabla$ on it. Roughly speaking it is a group of transforms of the tangent vectors generated by parallel transport along all closed loops.
The Ambrose-Singer theorem states that the Lie algebra $\mathfrak{hol}$ of Lie group $Hol(g)$ is generated by the Riemann curvature tensor $R$ and all it's derivatives at point $\gamma(0)$: one should consider an algebra generated by anti-symmetric matrices $\{R(X,Y),~ (\nabla R)(X,Y),\ldots ,(\nabla ^n R)(X,Y), \ldots\}$ for all vectors $X,Y \in T_{\gamma(0)}M$ and for all $n\geq 0$. The curvature tensor itself can be expressed through the connection:
$$
R(\xi,\eta)\zeta= \nabla_{\xi}\nabla_{\eta} \zeta - \nabla_{\eta}\nabla_{\xi} \zeta -\nabla_{ [ \xi,\eta ] } \zeta,
$$
this formula is a standard definition of Riemann tensor. Also in Riemannian geometry there is a 'second variation formula' that express the second variation of the energy functional $E$ on the geodesics $\gamma$ along vector fields $W_1, W_2$ via Riemann tensor:
$$
\frac{\partial^2 E(\sigma (u_1,u_2))}{\partial u_1 \partial u_2}\Big| _{u_1=u_2=0}=-2\int \langle W_2, \frac{D^2 W_1}{d t^2} + R(W_1,\dot{\gamma})\dot{\gamma} \rangle dt,
$$
where $\sigma=\sigma(t,u_1,u_2)$ is a two-parametric variation of the geodesic $\gamma(t)$, $W_i=\frac{\partial \sigma (t,u_1,u_2)}{\partial u_i}\Big |_{u_1=u_2=0}$ --- vector fields along $\gamma$ and $D$ is a restriction of $\nabla$ on $\gamma$.

From general point of view one can say that curvature $R$ (which is a second-order operator of $g$ and which determines almost all geometry of the space) can be defined via first-order operator $\nabla$. In some sense in Riemannian geometry we always can decrease the order of a problem using Levi-Civita connection. In subRiemannian geometry the situation is much more complicated. There is no privileged connection, at least at the present moment, and this fact gives us another reason to call subRiemannian geometry non-integrable.

Still in many cases one can define a curvature knowing some specific characteristics of the sub-Riemannian structure. One can analyze the behavior of normal geodesics for some concrete sub-Riemannian structures and using the 'second variation formula' extract from this data a definition of the curvature. This was done in \cite{AGL}.

In this paper we define the holonomy flag for sub-Riemannian structure using Ambrose-Singer theorem via formal curvature tensor. We present the analog of Codazzi and Weingarten equations for 3-dimensional Sasakian models and discuss some ideas how the covariant derivative should  for the .

The author is grateful to A.A. Agrachev, D.V. Alekseevsky and Ya.V. Bazaikin for the stimulating conversations.

\section{Codazzi equations.}

Let us remind some general facts from the differential geometry of 2-dimensional surfaces in Euclidian $\mathbb{R}^3$. The surface $\Sigma$ in $\mathbb{R}^3$ is the map $r:U\rightarrow \mathbb{R}^3$ where $U$ --- some open subset in $\mathbb{R}^2$. The vector fields $r_i=\frac{\partial r}{\partial u^i}$ for $i=1,2$ forms the basis of tangential plane $T_p \Sigma$ at every point $p\in \Sigma$. In contemporary differential geometry tangent vector fields $r_i$ usually identified with partial derivatives $\frac{\partial}{\partial u^i}$ along coordinates $u^i$ although usually there is no immersion $r$. Next we define the first quadratic form
$$
I=g_{ij}\Big| _{i,j=1,2}=\left(
  \begin{array}{cc}
    \langle r_1, r_1\rangle & \langle r_1, r_2\rangle \\
    \langle r_2, r_1\rangle & \langle r_2, r_2\rangle \\
  \end{array}
\right),
$$
which is just a scalar product induced in the tangent plane $T_p \Sigma$ by Euclidean product $\langle \cdot, \cdot \rangle$ in $\mathbb{R}^3$. This is a metric.
Then we define a unitary normal vector field $n=\frac{r_1 \times r_2}{|r_1 \times r_2|}$, where $\times$ is a vector product in $\mathbb{R}^3$. The second fundamental form  is
$$
II=b_{ij}\Big| _{i,j=1,2}=\left(
  \begin{array}{cc}
    \langle r_{11}, n\rangle & \langle r_{12}, n\rangle \\
    \langle r_{21}, n\rangle & \langle r_{22}, n\rangle \\
  \end{array}
\right),
$$
where $r_{ij}=\frac{\partial^2 r}{\partial u^i \partial u^j}$. The Gaussian curvature is $K=\frac{det ~II}{det ~I}$.

We have a frame $\{r_1,r_2,n\}$ which is a basis of $T_p\mathbb{R}^3=\mathbb{R}$ for $p\in \Sigma$. We emphasize that this frame is parametrised by points $p$ on $\Sigma$ and outside of $\Sigma$ does not exist. As $\{r_1,r_2,n\}$ is a basis  we can consider a decomposition
$$
r_{ij} = \Gamma_{ij}^k r_k + b_{ij}n. \eqno{(1)}
$$
This are the \emph{derivation equations}. The Christoffel symbols $\Gamma_{ij}^k$ are defined from (1) as coefficients in this decomposition. This symbols define the connection completely and can be calculated via coefficients $g_{ij}$. The \emph{Weingarten equations} describe how does the normal $n$ deforms while the frame moves along the surface:
$$
n_i= -b_{ij}g^{jk}r_k.
$$
The derivation equations and Weingarten can be gathered into two matrix systems:
$$
\frac{\partial}{\partial u^i}
\left(
  \begin{array}{c}
    r_1 \\
    r_2 \\
    n \\
  \end{array}
\right)=A_i \left(
  \begin{array}{c}
    r_1 \\
    r_2 \\
    n \\
  \end{array}
\right),
$$
where
$$
A_1=\left(
  \begin{array}{ccc}
    \Gamma_{11}^1 & \Gamma_{11}^2 & b_{11} \\
    \Gamma_{12}^1 & \Gamma_{12}^2 & b_{12} \\
    -b_{1j}g^{j1} & -b_{1j}g^{j2} & 0 \\
  \end{array}
\right),~~
A_2=\left(
  \begin{array}{ccc}
    \Gamma_{21}^1 & \Gamma_{21}^2 & b_{21} \\
    \Gamma_{22}^1 & \Gamma_{22}^2 & b_{22} \\
    -b_{2j}g^{j1} & -b_{2j}g^{j2} & 0 \\
  \end{array}
\right).
$$
We remind that $(u^1,u^2)$ are the standard coordinates on some domain $U$ in $\mathbb{R}^2$. So the partial derivatives along them will always commute. This leads to the following equations on the matrices $A_i$:
$$
\frac{\partial}{\partial u^2}A_1 -\frac{\partial}{\partial u^1}A_2 +A_1A_2-A_2A_1=0.\eqno{(2)}
$$
This are the \emph{Codazzi equations}. Equations (2) can be interpreted in a different way. Let us consider new derivatives $\nabla_i=\frac{\partial}{\partial u^i}-A_i$, then their commutativity is equivalent to (2). It is easy to see that the following quantity is a curvature:
$$
R=[\nabla_1,\nabla_2]=\nabla_1 \nabla_2 - \nabla_2 \nabla_1.
$$
$R$ is an operator acting on the frame $\{r_1,r_2,n\}$. It's clear that $R$ is identically equal to zero if the matrices $A_i$ were generated by some surface $r:U\rightarrow \mathbb{R}^3$. For the subRiemannian case we will show that $R$ can be identified with a skew-symmetric non-vanishing $3\times3$-matrix.

\vskip5mm

Next we want to generalize this construction to sub-Riemannian 3-dimensional cases. We consider frame $(X,Y,Z)$ in $\mathbb{R}^3$ (or $S^3$) instead of $(r_1,r_2,n)$. Where $(X,Y)$ is a bracket-generating $(2,3)$-distribution such that $[X,Y]=Z$. The main difference here is that the distribution $(X,Y)$ is non-integrable and there is no 2-dimensional surface $\Sigma$ with coordinates $(u^1,u^2)$. In classical case $T\Sigma$ forms the set of all admissible directions and we will identify $r_1$ with $X$ and $r_2$ with $Y$ though $X$ and $Y$ are not tangent to any surface. Next we should decide how will we generalize the derivative $\frac{\partial}{\partial u^i}r_j$. The easiest way is to consider a Lie derivative $\mathfrak{L}_XY$. We will illustrate that this approach can give some curious results, but we want to consider such a generalization that if $[X,Y]=0$ then our construction should become a classical one. So we can not use a Lie derivative because $\mathfrak{L}_XX=0$ (that should be the analog for $r_{11}$) for any $X$.

\section{Heisenberg and $SU(2)$.}

In this section we will carry out the previous construction for two subRiemannian structures: Heisenberg group and $SU(2)$. For Heisenberg
$$
X=\frac{\partial}{\partial x}-\frac{y}{2}\frac{\partial }{\partial z}=\left(
  \begin{array}{c}
    1 \\
    0 \\
    -\frac{y}{2} \\
  \end{array}
\right),\quad Y=\frac{\partial}{\partial y}+\frac{x}{2}\frac{\partial }{\partial z}=\left(
  \begin{array}{c}
    0 \\
    1 \\
    \frac{x}{2} \\
  \end{array}
\right), \quad Z=\frac{\partial}{\partial z}=\left(
  \begin{array}{c}
    0 \\
    0 \\
    1 \\
  \end{array}
\right).
$$
It is easy to check that $[X,Y]=Z$. Let us calculate the derivative equations firstly. By $\delta_U V$ we denote the standard derivative of vector field $V$ along the vector field $U$, this is the generalization of $\frac{\partial}{\partial u^j} r_i$.
$$
\delta _X X=(\frac{\partial}{\partial x}-\frac{y}{2}\frac{\partial }{\partial z})\left(
  \begin{array}{c}
    1 \\
    0 \\
    -\frac{y}{2} \\
  \end{array}
\right)=0,\quad \delta _X Y=(\frac{\partial}{\partial x}-\frac{y}{2}\frac{\partial }{\partial z})\left(
  \begin{array}{c}
    0 \\
    1 \\
    \frac{x}{2} \\
  \end{array}
\right)=\left(
  \begin{array}{c}
    0 \\
    0 \\
    \frac{1}{2} \\
  \end{array}
\right),
$$
$$
\delta _Y X=(\frac{\partial}{\partial y}+\frac{x}{2}\frac{\partial }{\partial z})\left(
  \begin{array}{c}
    1 \\
    0 \\
    -\frac{y}{2} \\
  \end{array}
\right)=\left(
  \begin{array}{c}
    0 \\
    0 \\
    -\frac{1}{2} \\
  \end{array}
\right),\quad \delta _Y Y=(\frac{\partial}{\partial y}+\frac{x}{2}\frac{\partial }{\partial z})\left(
  \begin{array}{c}
    0 \\
    1 \\
    \frac{x}{2} \\
  \end{array}
\right)=0,
$$
The Weingarten equations take the form:
$$
\delta_X Z=\delta_Y Z =0.
$$
Then
$$
A_1=\left(
  \begin{array}{ccc}
    0 & 0 & 0 \\
    0 & 0 & \frac{1}{2} \\
    0 & 0 & 0 \\
  \end{array}
\right),~~
A_2=\left(
  \begin{array}{ccc}
    0 & 0 & -\frac{1}{2} \\
    0 & 0 & 0 \\
    0 & 0 & 0 \\
  \end{array}
\right).
$$
And it easy to check that the curvature $R$ is identically vanish for the Heisenberg group, it confirms that this geometry is the flatest among the 3D subRiemannian geometries. We also note that all the Christoffel symbols are vanish, for the Riemannian geometry it would automatically mean that the space is flat. Here we have that the components $b_{ij}$ do not form the quadratic form, they appear only as a coefficients of decomposition (1).

\vskip3mm

Next we will derive the same formulas for the $S^3=SU(2)$ assuming that $S^3=\{ q= x+iy+jz+kt\in \mathbb{H}\big| |q|=1\}$. We will use the following coordinates on $S^3$
$$
\left(
  \begin{array}{c}
    x \\
    y \\
    z \\
    t \\
  \end{array}
\right) =
\left(
  \begin{array}{c}
    \cos(\theta) \cos(\psi +\phi)\\
    \cos(\theta) \sin(\psi +\phi)\\
    \sin(\theta) \cos(\psi -\phi)\\
    \sin(\theta) \sin(\psi -\phi)\\
  \end{array}
\right).
$$
The vectors $\{iq,jq,kq\}$ will be orthogonal to the vector $q\in S^3$ and will be a orthonormal basis of $T_q S^3$. We will express them through the basis of coordinates vector fiels:
$$
\partial_{\phi}\equiv \frac{\partial q}{\partial \phi} =\left(
  \begin{array}{c}
    -\cos(\theta) \sin(\psi +\phi)\\
    \cos(\theta) \cos(\psi +\phi)\\
    \sin(\theta) \sin(\psi -\phi)\\
    -\sin(\theta) \cos(\psi -\phi)\\
  \end{array}
\right),
\partial_{\theta}\equiv \frac{\partial q}{\partial \theta} =\left(
  \begin{array}{c}
    -\sin(\theta) \cos(\psi +\phi)\\
    -\sin(\theta) \sin(\psi +\phi)\\
    \cos(\theta) \cos(\psi -\phi)\\
    \cos(\theta) \sin(\psi -\phi)\\
  \end{array}
\right),
$$
$$
\partial_{\psi}\equiv \frac{\partial q}{\partial \psi} =\left(
  \begin{array}{c}
    -\cos(\theta) \sin(\psi +\phi)\\
    \cos(\theta) \cos(\psi +\phi)\\
    -\sin(\theta) \sin(\psi -\phi)\\
    \sin(\theta) \cos(\psi -\phi)\\
  \end{array}
\right).
$$
Note that the frame $\{\partial_{\phi},\partial_{\theta},\partial_{\psi}\}$ is not orthonormal. Then we denote the vector $iq=-y+ix-jt+kz$ as $Z$ and we see that $Z=\partial_{\psi}=\left( \begin{array}{c} 0 \\0 \\1\end{array}\right)$. After solving two linear system one can get
$$
X=kq=-\frac{\cos(2\psi)}{\sin(2\theta)}\partial_{\phi} + \sin(2\psi)\partial_{\theta} +\frac{\cos(2\theta)\cos(2\psi)}{\sin(2\theta)}\partial_{\psi}=\left(
  \begin{array}{c}
   -\frac{\cos(2\psi)}{ \sin(2\theta)}\\
   \sin(2\psi)\\
   \frac{\cos(2\theta)\cos(2\psi)}{\sin(2\theta)}\\
  \end{array}
\right),
$$
$$
Y=jq=\frac{\sin(2\psi)}{\sin(2\theta)}\partial_{\phi} + \cos(2\psi)\partial_{\theta} -\frac{\cos(2\theta)\sin(2\psi)}{\sin(2\theta)}\partial_{\psi}=\left(
  \begin{array}{c}
  \frac{\sin(2\psi)}{\sin(2\theta)}\\
  \cos(2\psi)\\
  -\frac{\cos(2\theta)\sin(2\psi)}{\sin(2\theta)}\\
  \end{array}
\right).
$$

The first derivative equation is:
$$
\delta_X X=-\frac{\cos(2\psi)}{\sin(2\theta)}\partial_{\phi} \left(
  \begin{array}{c}
   -\frac{\cos(2\psi)}{ \sin(2\theta)}\\
   \sin(2\psi)\\
   \frac{\cos(2\theta)\cos(2\psi)}{\sin(2\theta)}\\
  \end{array}
\right)+ \sin(2\psi)\partial_{\theta}\left(
  \begin{array}{c}
   -\frac{\cos(2\psi)}{ \sin(2\theta)}\\
   \sin(2\psi)\\
   \frac{\cos(2\theta)\cos(2\psi)}{\sin(2\theta)}\\
  \end{array}
\right)
$$
$$
+\frac{\cos(2\theta)\cos(2\psi)}{\sin(2\theta)}\partial_{\psi}\left(
  \begin{array}{c}
   -\frac{\cos(2\psi)}{ \sin(2\theta)}\\
   \sin(2\psi)\\
   \frac{\cos(2\theta)\cos(2\psi)}{\sin(2\theta)}\\
  \end{array}
\right)=\left(
  \begin{array}{c}
   4\frac{\cos(2\psi)\sin(2\psi)\cos(2\theta)}{ \sin(2\theta)^2}\\
   2\frac{\cos(2\psi)^2\cos(2\theta)}{\sin(2\theta)}\\
  -2\frac{\cos(2\psi)\sin(2\psi)(1+\cos(2\theta)^2)}{\sin(2\theta)^2}\\
  \end{array}
\right)=
$$
$$
=\Gamma_{11}^1X +\Gamma_{11}2 Y +b_{11}Z,
$$
where $\Gamma_{11}^1=-2\frac{\cos(2\psi)^2 \sin(2\psi)\cos(2\theta)}{\sin(2\theta)}, ~\Gamma_{11}^2=2\frac{\cos(2\psi)\cos(2\theta)(1+\sin(2\psi)^2)}{\sin(2\theta)}$ and $b_{11}=-2\cos(2\psi)\sin(2\psi)$.

Making the same calculations for $\delta_X Y,~\delta_Y X$ and $\delta_Y Y$ we get the following quantities:
$$
\Gamma_{12}^1=\Gamma_{21}^1=-2\frac{\cos(2\theta)\cos(2\psi)^3}{\sin(2\theta)},~
\Gamma_{12}^2=\Gamma_{21}^2=-2\frac{\cos(2\theta)\sin(2\psi)^3}{\sin(2\theta)},
$$
$$
\Gamma_{22}^1=2\frac{\cos(2\theta)\sin(2\psi)(1+\cos(2\psi)^2)}{\sin(2\theta)},~
\Gamma_{22}^2=-2\frac{\cos(2\theta)\sin(2\psi)^2\cos(2\psi)}{\sin(2\theta)}.
$$
The components $b_{ij}$ are:
$$
~b_{12}=2\sin(2\psi)^2,~ b_{21}=-2\cos(2\psi)^2,~b_{22}=2\sin(2\psi)\cos(2\psi).
$$
Note that here (as in the Riemannian case) $\Gamma_{12}^k=\Gamma_{21}^k$, but $b_{12}\neq b_{21}$ and it is easy to see that $[X,Y]=\delta_XY-\delta_YX = b_{12}Z-b_{21}Z=2Z$.
Weingarten equations are trivial $\delta_XZ=\delta_YZ=0$. Using this coefficients one can form matrices $A_1$ and $A_2$ and check that
$$
R=\delta_Y A_1 -\delta_X A_2 +A_1A_2-A_2A_1=\left(
                                               \begin{array}{ccc}
                                                 0 & -4 & 0 \\
                                                 4 & 0 & 0 \\
                                                 0 & 0 & 0 \\
                                               \end{array}
                                             \right).
$$
This result is not surprising because this $R$ corresponds to a sectional curvature of two-dimensional base $S^2$ in the Hopf fibration. Still we can consider $R$ as a curvature of the subRiemannian structure generated by non-holonomic distribution $\{X,Y,Z\}$. And as wee see $R$ is zero for the Heisenberg and it is constant for $SU(2)$.

\vskip3mm

Now let us make the same calculations using the Cartan method, for this one should consider the Lie derivative $\mathfrak{L}$ instead of standard derivative $\delta$. Consider the frame $\{X,Y,Z\}$ such that $[X,Y]=2Z,~[Y,Z]=\rho X,~[X,Z]=-\rho Y$, where $\rho\in\{0,-1,+1\}$. Then the derivative equations are
$$
\mathfrak{L}_X X=[X,X]=0,~ \mathfrak{L}_X Y=2Z,~ \mathfrak{L}_Y X=-2Z, ~\mathfrak{L}_YY=0.
$$
And the Weingarten equations are
$$
\mathfrak{L}_X Z=-\rho Y,~~ \mathfrak{L}_Y Z=\rho X.
$$
The matrices:
$$
A_1=\left(
  \begin{array}{ccc}
    0 & 0 & 0 \\
    0 & 0 & 2 \\
    0 & -\rho & 0 \\
  \end{array}
\right),~~
A_2=\left(
  \begin{array}{ccc}
    0 & 0 & -2 \\
    0 & 0 & 0 \\
    \rho & 0 & 0 \\
  \end{array}
\right).
$$
The curvature $R$ is
$$
R=\mathfrak{L}_Y A_1 -\mathfrak{L}_X A_2 +A_1A_2-A_2A_1=\left(
                                               \begin{array}{ccc}
                                                 0 & -2\rho & 0 \\
                                                 2\rho & 0 & 0 \\
                                                 0 & 0 & 0 \\
                                               \end{array}
                                             \right).
$$
And it coincides with the one for Heisenberg ($\rho=0$) and the one for $SU(2)$ ($\rho=2$).

\vskip3mm

Let us remind that in Riemannian geometry there is a formula that allow to calculate Christoffel symbols from (1) using the first quadratic form:
$$
\Gamma_{ij}^k=\frac{1}{2}g^{kl}\Big(\frac{\partial}{\partial u^i} g_{jl} +\frac{\partial}{\partial u^j} g_{il} -\frac{\partial}{\partial u^l} g_{ij}\Big).\eqno{(3)}
$$
Our admissible vector fields were orthonormal, the first quadratic form is identity $2\times 2$ -matrix and all the Christoffel symbols will vanish if one woild like to use (3) or some variant of this formula. Also for the Cartan method the Christoffel symbols are zero.

The main purpose of the calculations above was to illustrate that the Christoffel symbols could not be calculated in subRiemannian case without using the coordinate form of the frame $\{X,Y,Z\}$ and could not be calculated using the metric on the horizontal subspace generated by admissible directions $X$ and $Y$.

So up until now we are very pessimistic about defining some canonical and invariant connection for subRiemannian geometries without using some principal bundle of frames (see, for example, \cite{Vershik}).

\section{Holonomy algebra for non-holonomic distributions.}

In this section we suggest a new definition of the holonomy algebra for subRiemannian structures in 3D case. The usual way to define the holonomy group of Riemannian manifold uses the Levi-Civita connection $\nabla$ on it. This connection is unique and depend only on a metric $g$. If there is a loop $\gamma(t),~ t\in[0,1]$ on a manifold $M$ one should consider a tangent vector $v\in T_{\gamma(0)}M$ and parallelly transport it along the loop $\gamma$ using Levi-Civita connection. After this transport one will get tangent vector $\phi_{\gamma}(v)\in T_{\gamma(1)}M$. This map $\phi_{\gamma}$ will be an element in $SO(n)$ because Levi-Civita connection preserves the lengths and angles along the curves. If one will consider all loops which begin and end in point $\gamma(0)$ then he will get a closed subgroup $Hol(g)$ in $SO(n)$, it is easy to prove that this subgroup does not depend (up to a conjugation) on the choice of fixed point $\gamma(0)=\gamma(1)$ if $M$ was connected. The Ambrose-Singer theorem states that the Lie algebra $\mathfrak{hol}$ of Lie group $Hol(g)$ is generated by the Riemann curvature tensor $R$ and all it's derivatives at point $\gamma(0)$: one should consider an algebra generated by anti-symmetric matrices $\{R(X,Y),~ (\nabla R)(X,Y),\ldots ,(\nabla ^n R)(X,Y), \ldots\}$ for all vectors $X,Y \in T_{\gamma(0)}M$ and for all $n\geq 0$. Obviously, the more flat the manifold --- the smaller holonomy group it has. The curvature tensor itself can be expressed through the connection:
$$
R(X,Y)Z= \nabla_{X}\nabla_{Y} Z - \nabla_{Y}\nabla_{X} Z -\nabla_{ [ X,Y ] } Z, \eqno{(4)}
$$
this formula is a standard definition of Riemann tensor. The formula (4) has very simple geometric interpretation: vector $R(X,Y)Z$ is a result of a parallel transport of the vector $Z$ along infinitely small parallelogram with sides $X,Y,-X$ and $-Y$.

For the left-invariant metric on the Lie group $G$ and left-invariant vector fields $X,Y$ and $Z$ the following formula is true (\cite{Besse}, 7.21 and 7.28)
$$
\nabla_X Y= \frac{1}{2}[X,Y]  + U(X,Y), \eqno{(5)}
$$
where
$$
2\langle U(X,Y), Z \rangle =\langle[Z,X],Y\rangle + \langle X,[Z,Y]\rangle \eqno{(6)}
$$
and $\langle\cdot , \cdot \rangle$ is a left-invariant scalar product on the Lie algebra $\mathfrak{g}$. One can notice the resemblance between the formula (5) and the standard definition of covariant derivative of the vector field $V$ along the vector field $\frac{\partial}{\partial x_i}=e_i$:
$$
\nabla_{e_i} V=\Big(\frac{\partial V^k}{\partial x_i} + \Gamma^k_{ji}V^j\Big)e_k,
$$
where the partial derivative can be identified with the Lie bracket and the Christoffel symbols can be identified with affine shift $U$.



For the subRiemannian manifold $M$ one has the \emph{set of admissible directions} $\mathcal{H}=\mathcal{H}^0$. Usually it is a non-integrable subbundle of the tangent bundle $TM$ and called the horizontal distribution. Then $\mathcal{H}_1=[\mathcal{H},\mathcal{H}]$, $\mathcal{H}_2=[\mathcal{H},\mathcal{H}_1]+\mathcal{H}_1, \dots ,\mathcal{H}_i=[\mathcal{H},\mathcal{H}_{i-1}]+\mathcal{H}_{i-1}$. Finally, for some $k$ one has $\mathcal{H}_k=TM$. Obviously $\mathcal{H}^i\subset\mathcal{H}^{i+1}$.

We want to define the curvature for the subRiemannian structures also as commutation of the covariant derivatives (4). But we will show that on Lie algebra $\mathfrak{g}$ there is no such scalar product that projectives to scalar product on the distribution and generates the correct connection.

Now we construct the object which generalizes the holonomy algebra for the subRiemannian case. We called it \emph{holonomy flag}. Let $X,~Y\in \mathcal{H}$ be two horizontal vector fields. Consider the family of subspaces $\mathcal{R}^i_{XY}=\{R(X,Y)Z~|~Z\in \mathcal{H}^i\}\bigcap \mathcal{H}^i$ in $TM$. Again, obviously $\mathcal{R}^i\subset \mathcal{R}^{i+1}$.

To make this object more concrete we fix the frames $X^i_j$ of the $\mathcal{H}^i$ and $\mathcal{R}^i_{XY}=span\{ R(X,Y)X^i_j| j=0\ldots dim(\mathcal{H}^i)\}$. The collection of all this chains of the vector subspaces for all pairs of admissible vector fields $\mathcal{R}=\{\mathcal{R}^i_{X_jX_k}|X_j,X_k\in \mathcal{H}, j<k=1\ldots dim(\mathcal{H})\}$ forms the holonomy flag. We emphasize that we consider only only admissible vector fields $X,Y$ in the definition of $\mathcal{R}_{XY}$ because they play role of the sides of a parallelogram along which the field $Z$ is transported.

In the Riemannian case there is invertible metric $G$ and one can consider $R(X,Y)$ as an anti-symmetric map acting on the tangent space. But in the subRiemannian case $G^{-1}$ is degenerate and it is impossible to lift the indices in $R$ up or to put them down. For the Riemannian manifold, as $\mathcal{H}^0=TM$, one get $\mathcal{R}=\mathcal{R}^0=\{R(X_i,X_j)X_k|i,j,k=1\ldots dim(M)\}$, by multiplying this set of 3-linear maps by $G^{-1}$ one will get the set of skew-symmetric matrices, which will generate the holonomy algebra $\mathfrak{hol}$. It means that this construction is natural and coincides with Riemannian one in this trivial case.




For any subRiemannian structure on 3D Lie group one can chose frame such that the Lie brackets will take the following form (\cite{Agrachev-Barilari})
$$
[X,Y]=Z,\quad [Y,Z]=(\chi + \kappa)X,\quad [X,Z]=(\chi-\kappa)Y. \eqno{(7)}
$$
For this structure $\mathcal{H}^0=span\{X,Y\}$, and we declare the vector fields $X$ and $Y$ to be orthonormal with respect to the scalar product $\langle\cdot ,\cdot \rangle_0$. We define the connection by the formula (5), where the Lie bracket is defined by (7). Symmetric and bilinear $U$ we will define as follows:
$$
U(X,Z)=\alpha Z, \quad U(Y,Z)=\beta Z
$$
for some fixed constants $\alpha, \beta \in \mathbb{R}$. Then the connection will take form
$$
\nabla_X Y=\frac{1}{2} Z, \quad \nabla_Y X=-\frac{1}{2}Z,
$$
$$
\nabla_X Z=\frac{\chi-\kappa}{2}Y + \alpha Z, \quad \nabla_Z X= \frac{\kappa-\chi}{2}Y +\alpha Z,\eqno{(8)}
$$
$$
\nabla_Y Z=\frac{\chi+\kappa}{2}X + \beta Z, \quad \nabla_Z Y= -\frac{\kappa+\chi}{2}X +\beta Z.
$$

Firstly we prove the following

\textbf{Lemma.} \emph{For the connection (8) the metric $\langle\cdot ,\cdot \rangle_0$ on the horizontal distribution is parallel with respect to admissible directions and has no torsion. There is no metric $\langle\cdot ,\cdot \rangle$ on the whole Lie algebra $\mathfrak{g}$ which projectivise to the metric $\langle\cdot ,\cdot \rangle_0$ and generates the $U$-term of connection (8) by condition (6).}

\textbf{Proof}  The torsionless is straight-forward. The torsion tensor is vanish because it is calculated by
$$
T(X,Y)=\nabla_X Y -\nabla_Y X -[X,Y],
$$
and $U$ is symmetric. The parallelness means that $\nabla \langle \cdot ,\cdot \rangle_0=0$ or $\langle \nabla_A B,C \rangle_0 + \langle B,\nabla_A C \rangle_0=0$ for any $A,B,C\in \mathcal{H}^0$. The latest equality is trivial consequence from (8) and the fact that $X$ and $Y$ are orthonormal. Also here we should demand for the scalar product $\langle \cdot, \cdot \rangle_0$ to be defined \emph{only} on $\mathcal{H}^0$ and formally $\langle Z,X\rangle_0=\langle Z,Y\rangle_0=0$.

Suppose that there is a metric $\langle \cdot, \cdot \rangle$ on $\mathcal{H}^1=TM$ generating the $U$-term. But then $0=2\langle U(X,Y),Z\rangle = (\kappa-\chi)\langle Y,Y\rangle - (\kappa+\chi)\langle X,X \rangle=-2\chi$ because we require for the metric $\langle \cdot,\cdot \rangle$ to projectivise to the metric $\langle \cdot, \cdot\rangle_0$ and $\langle X,X\rangle=\langle Y,Y \rangle=1$. So it is possible only for the case when $\chi=0$. Lemma proved.

\vskip5mm

This Lemma shows that the defined connection is natural from the standard point of view in Riemannian geometry that the connection should be torsionless and the metric should be parallel with respect to this connection. We emphasize here that we could define $U$-term in a different way, for example, we could put $U(X,Y)=-\chi Z$ and then the calculation in the Lemma could be satisfied, but in this case the curvature will be spoiled. Next we prove that the curvature (and the holonomy flag, respectively) are generalized in an appropriate way too.

\vskip5mm

\textbf{Theorem.} \emph{Let $M$ be the subRiemannian 3D-Lie group with frame $\{X,Y,Z\}$ and the structure of Lie brackets (6). Then the holonomy flag defined by connection (7) is
$$
\mathcal{R}^0=span\{\frac{1}{4}(\chi-\kappa)Y, \frac{1}{4}(\kappa +\chi)X\},
$$
$$
\mathcal{R}^1=span\{\frac{1}{4}(\chi-\kappa)Y-\frac{3}{2}\alpha Z, \frac{1}{4}(\kappa +\chi)X-\frac{3}{2}\beta Z,-\frac{\alpha}{2}(\chi+\kappa)X+\frac{\beta}{2}(\chi-\kappa)Y\}.
$$
}
\textbf{Proof}

The proof is straight-forward.
$$
R(X,Y)X=\nabla_X \nabla_Y X -\nabla_Y \nabla_X X-\nabla_{[X,Y]} X=-\frac{1}{2}\nabla_X Z -\nabla_Z X=
$$
$$
=-\frac{1}{2}(\frac{\chi-\kappa}{2}Y+\alpha Z) - \frac{\kappa-\chi}{2}Y - \alpha Z=\frac{\chi-\kappa}{4}Y -\frac{3}{2}\alpha Z,
$$
$$
R(X,Y)Y=\nabla_X \nabla_Y Y -\nabla_Y \nabla_X Y-\nabla_{[X,Y]} Y=-\frac{1}{2}\nabla_Y Z -\nabla_Z Y=
$$
$$
=-\frac{1}{2}(\frac{\chi+\kappa}{2}X + \beta Z) - \frac{\kappa+\chi}{2}X -\beta Z=\frac{\chi+\kappa}{4}X -\frac{3}{2}\beta Z,
$$
and
$$
R(X,Y)Z=\nabla_X \nabla_Y Z -\nabla_Y \nabla_X Z-\nabla_{[X,Y]} Z=\nabla_X (\frac{\chi+\kappa}{2}X + \beta Z)
$$
$$
-\nabla_Y (\frac{\chi-\kappa}{2}Y + \alpha Z)=\beta\nabla_X Z- \alpha\nabla_Y Z= \beta(\frac{\chi-\kappa}{2}Y + \alpha Z)-\alpha(\frac{\chi+\kappa}{2}X + \beta Z)=
$$
$$
=-\frac{\alpha}{2}(\chi+\kappa)X+\frac{\beta}{2}(\chi-\kappa)Y.
$$
Then by considering the projective of $R(X,Y)X$ and $R(X,Y)Y$ on $\mathcal{H}^0$ one will get $\mathcal{R}^0$. Theorem proved.

\vskip5mm

One can consider the tangent space $T_p M $ as a flat approximation of the Riemannian manifold $M$ at point $p$. It is known that the Heisenberg group is a nilpotent approximation for any subRiemannian 3D Lie group \cite{Agrachev-Barilari}. We see that for the Heisenberg group, i.e. for $\chi=\kappa=0$, the holonomy flag vanishes; just like the curvature vanishes for the tangent space $T_p M$. That give us an evidence that our definition of the holonomy flag is quite natural.

We conclude with the following observation. If one will put $\alpha=\beta=0$ then the connection (8) will coincide with connection on 3D Lie groups with bi-invariant metric. It is the most privileged metric on Lie groups. Unfortunately this metric is well-defined only for the compact groups. In subRiemannian case we see that appropriate connection stays reasonable.

\end{document}